\documentclass[a4paper,12pt]{article}
\usepackage{amssymb,amsmath,latexsym}

\usepackage{wrapfig}
\usepackage[dvips]{graphicx}
\usepackage{amsmath,amsthm}
\usepackage{amsfonts,amssymb}
\usepackage{amscd}
\usepackage{array}
\usepackage{amssymb,latexsym}
\usepackage {amsfonts, amsmath, amssymb, amsthm, color, wasysym}
\usepackage[matrix, arrow, curve]{xy}

\def\mat{\begin{pmatrix} }
\def\endm{\end{pmatrix} }

\newtheorem{theorem}{Theorem}[section]
\newtheorem{lemma}[theorem]{Lemma}
\newtheorem{proposition}[theorem]{Proposition}

\theoremstyle{definition}

\newtheorem{example}[theorem]{Example}
\newtheorem{definition}[theorem]{Definition}
\newtheorem{remark}[theorem]{Remark}

\newtheorem{corollary}[theorem]{Corollary}

\DeclareMathOperator{\Soc}{Soc} 

\DeclareMathOperator{\Next}{Next}

 \DeclareMathOperator{\Sb}{SB}
  \DeclareMathOperator{\Cw}{CW}

\def\lan{\langle}
\def\ran{\rangle}

\def\lan{\langle}
\def\ran{\rangle}

\begin{document} 


\author{Mikhail Antipov
\\[-1pt]
\small Saint-Petersburg University\\[-3pt]
{\tt \small hyperbor@list.ru}
}
\title{Derived equivalence of symmetric special biserial algebras}

\maketitle{}

\begin{abstract}
We introduce Brauer complex of symmetric SB-algebra, and reformulate in terms of Brauer complex the so far known
invariants of stable and derived equivalence of symmetric SB-algebras. In particular, the genus
of Brauer complex turns out to be invariant under derived equivalence. We study transformations
of Brauer complexes which preserve class of derived equivalence. Additionally, we establish a new invariant of derived
equivalence of symmetric SB-algebras. As a consequence, symmetric SB-algebras with Brauer complex of genus 0 are
classified.

\end{abstract}

\tableofcontents

\section{Introduction}

The present paper lies within a series of papers, devoted to classification of symmetric special biserial
algebras up to derived equivalence (i.e., up to equivalence of derived categories).
Recall that a symmetric SB-algebra $\Lambda$ is uniquely determined by a pair $(\Gamma(\Lambda), f)$, where $\Gamma(\Lambda)$ is the Brauer graph of $\Lambda$ and $f:V(\Gamma(\Lambda))\to \mathbb{N}$ maps vertices of $\Gamma(\Lambda)$ to their multiplicities (see, e.g.,~\cite{1} and Proposition~\ref{bijection}).
\begin{itemize}
  \item We show that the multiset of multiplicities of vertices of $\Gamma(\Lambda)$ is invariant under derived equivalence (Proposition~\ref{multiplicities}). In order to prove this, we determine the center $Z(\Lambda)$.
  \item In section 3 we introduce {\it Brauer} $\Cw$-{\it complex} $C(\Lambda)$ --- a relevant tool for studying derived
  equivalence. Topologically, $C(\Lambda)$ is a sphere with handles. We reformulate in terms of $C(\Lambda)$ the basic
  notions related to algebra $\Lambda$ and the invariants of stable equivalence, which appeared in~\cite{2}.
  In particular, the genus of $C(\Lambda)$ turns out to be invariant under stable equivalence. By a celebrated
  theorem of Rickard \cite{6}, these invariants are invariants of derived equivalence, too.
  \item We introduce {\it elementary tilting complexes} over symmetric special biserial algebras --- a generalization 
  of tilting complexes, which were treated in~\cite{3} (section 4). Equivalences of algebras, corresponding to
  elementary tilting complexes, can be reformulated in terms of 'elementary transformations' of Brauer $\Cw$-complexes
  of these algebras (Proposition~\ref{correspondence}). One sees that the algebra, which corresponds to the $\Cw$-complex obtained from $C(\Lambda)$ by an
  elementary transformation, is derived equivalent to $\Lambda$. Thus we obtain a direct graphic way of proving derived
  equivalence.
  \item In the last section we show that if the geometric realization of $C(\Lambda)$ is a sphere,
  then the invariants which we discuss in this paper determine $\Lambda$ up to derived equivalence.
\end{itemize}

\section{The center $Z(\Lambda)$ and the multiplicities of $A$-cycles}
Let $\Lambda$ be a symmetric SB-algebra over field $K$. Consider an extended quiver $Q_e=Q_e(\Lambda)$. Consider the
partitions of its arrow set into $A$-cycles and into $G$-cycles
(see\cite{1}). Recall that $A$-cycles (and their multiplicities) correspond to the vertices of Brauer graph $\Gamma(\Lambda)$. We denote $A$-cycles by lower-case latine letters and denote vertices of $\Gamma(\Lambda)$ by the correspondent upper-case latine letters. 

Let $\{c_1,c_2,\dots, c_k\}$ be the set of $A$-cycles. For each $i=1,\dots, k$
consider a cyclic sequence $(\alpha_{i,1},\alpha_{i,2},\dots,\alpha_{i, l_i})$
of arrows of the cycle $c_i$. Let $f(c_1),f(c_2),\dots,f(c_k)\in \mathbb{N}$ denote the multiplicities of $A$-cycles.
For each loop $\alpha=\alpha_{i,k}$ which is not formal, set
$$
q_{\alpha}=(\alpha_{i,k+1}\alpha_{i,k+2}\dots,
\alpha_{i,l_i}\dots\alpha_{i,k})^{f(c_i)-1}
\alpha_{i,k+1}\alpha_{i,k+2}\dots,
\alpha_{i,l_i}\dots\alpha_{i,k-1}.
$$

\begin{proposition} \label{multiplicities}
{\bf 1.} The center $Z(\Lambda)$ is generated as a vector space over $K$ by $1$ and by the elements of the following three forms:

\begin{enumerate}
    \item[a.] Elements $ m_{i,t}=(\alpha_{i,1}\alpha_{i,2}\dots,\alpha_{i,l_i})^{t}
+(\alpha_{i,2}\alpha_{i,3}\dots,\alpha_{i,1})^{t}+\dots
+(\alpha_{i,l_i}\alpha_{i,1}\dots,\alpha_{i,l_i-1})^{t}$ for all
$i=1,2,\dots, k$ and $t=1,\dots , f(c_i)-1$.
\item[b.] Elements $q_{\alpha}$ for each non-formal loop $\alpha$.
  \item[c.] Elements $s_r=(\alpha_{i_r,1}\alpha_{i_r,2}\dots,\alpha_{i_r,l_{i_r}})^{f(c_{i_r})}$
for each vertex $r$ of $Q_e$, where $c_{i_r}$ is one of the two $A$-cycles, passing through $r$.
\end{enumerate}

  \smallskip\noindent {\bf 2.} $Z/(\Soc Z)\cong K[x_1,x_2,\dots,x_k]/\lan
\{x_i^{f(c_i)},(x_ix_j)_{i\neq j}\}
\ran$, where  i, j $\in$ {1, \dots, k}.

\smallskip \noindent {\bf 3.} The multiset $(f(c_1),f(c_2),\dots,f(c_k))$ is invariant under derived equivalence.
\end{proposition}
\begin{proof}
{\bf 1.} Recall that the value of $s_r$ doesn't depend on the choice of an $A$-cycle $c_{i_r}$ and that the elements $s_1, s_2,\dots ,s_n$ form a $K$-basis of $\Soc(\Lambda)$ (see, e.g.,~\cite{1}).
Since $\Lambda$ is a symmetric algebra, the socle $\Soc(\Lambda)$ is contained in $Z$, so $s_r\in Z$. Moreover, for a
non-formal loop $\alpha$ at vertex $r$ and for the corresponding idempotent $e_r$ and path $p \notin \{e_r, \mbox{ } \alpha\}$ we get
$e_rq_{\alpha}=q_{\alpha}=q_{\alpha}e_r$, $\alpha
q_{\alpha}=s_r=q_{\alpha}\alpha$, и $q_{\alpha}p=0=pq_{\alpha}$.
Thus $q_{\alpha}\in Z$.
Similarily, for all $i,t,r$ we get $e_rm_{i,t}=m_{i,t}e_r$, since the summands in $m_{i,t}$ are circuits. Furthermore, for all $l_1,l_2,t_1$
\begin{multline*}  
(\alpha_{i,l_1}\alpha_{i,l_1+1}\dots,\alpha_{i,l_1-1})^{t_1}\alpha_{i,l_1}\alpha_{i,l_1+1}\dots,\alpha_{i,l_2}m_{i,t}=\\
(\alpha_{i,l_1}\alpha_{i,l_1+1}\dots,\alpha_{i,l_1-1})^{t+t_1}\alpha_{i,l_1}\alpha_{i,l_1+1}\dots,\alpha_{i,l_2}= \\
m_{i,t}(\alpha_{i,l_1}\alpha_{i,l_1+1}\dots,\alpha_{i,l_1-1})^{t_1}\alpha_{i,l_1}\alpha_{i,l_1+1}\dots,\alpha_{i,l_2}\end{multline*}
\noindent
Since for the rest paths $p$ (subpaths of other $A$-cycles) $m_{i,t}p=pm_{i,t}=0$, we get $m_{i,t}\in Z$.

Each $z \in Z$ can be uniquely represented as
\begin{equation} \label{2}
z=\sum_{j=1}^{N} a_jp_j+s,
\end{equation}
\noindent where $0\neq a_j\in K$, paths $p_j$ are distinct nonzero paths in the quiver $Q_e$ which are not contained in the socle, $s\in \Soc (\Lambda)$.

By induction on the number of summands in the sum ~(\ref{2}) we show, that $z$ can be represented as a linear combination
of elements $m_{i,t}$ and $q_{\alpha}$.
Fix $i \in \{1, \dots N\}$ and write
\begin{equation*} 
p_i=\alpha_1\alpha_2\dots\alpha_{m},
\end{equation*}
\noindent where $\alpha_1,\alpha_2,\dots,\alpha_{m}$ are consequent arrows of an $A$-cycle $c_i$.
Let $\alpha_{m+1}$ be the next arrow of $c_i$.
There are two cases:

\smallskip
Case 1: $\alpha_1\alpha_2\dots\alpha_{m}\alpha_{m+1}\notin \Soc(\Lambda)$.
In this case the path $\alpha_{m+1}\alpha_1\alpha_2\dots\alpha_{m}$
has coefficient $a_i$ in the sum $\sum a_j\alpha_{m+1}p_j$. Since $z\alpha_{m+1}=\alpha_{m+1}z$, we obtain $\alpha_{m+1}=\alpha_1$, i.e.
$p_i=(\alpha_{u,1}\alpha_{u,2}\dots,\alpha_{u,l_u})^{t}$ for some $A$-cycle $c_u$, $t < f(c_u)$.
Moreover, the other summands of $m_{u,t}$ also have coefficient $a_i$ in the sum ~(\ref{2}).
We see that the sum representing element $z-a_im_{u,t}\in Z$ has less summands than the sum representing $z$,
so the inductive hypothesis is applied.

\smallskip
Case 2: $\alpha_1\alpha_2\dots\alpha_m\alpha_{m+1}=s_l\in \Soc(\Lambda)$ for some $l$.
Consider an idempotent $e_r$ such that $\alpha_m e_r\alpha_{m+1}\neq 0$.
The expressions for $ze_r$ and $e_rz$ must contain $p_i$ as a summand.
Therefore $p_i$ is a closed path. It follows that $\alpha_{m+1}$ is a loop and $p_i=q_{\alpha_{m+1}}$, and we apply the inductive hypothesis to $z-a_ip_i$.

\smallskip
{\bf 2.} Observe that $\Soc(Z)$ is generated by the elements $s_r$ and $q_{\alpha}$, for all loops $\alpha$ which are not
separate $A$-cycles (i.e., $\alpha q_{\alpha }=0$). Moreover,
$m_{i,t}m_{j,t_1}=\delta_{ij}m_{i,t+t_1}$ and $m_{i,t}^{f(c_i)}\in \Soc(Z)$. These two observations imply the claim.

\smallskip
{\bf 3.} The claim follows directly from p.{\bf 2}, since $Z(\Lambda)$ is invariant under derived equivalence
(see~\cite{5}). The maximal element $f(c_i)$ equals the maximal index of nilpotency of nilpotents in $Z/\Soc(Z)$; the remaining proof is by induction.

\end{proof}

\section{Brauer complex}
\subsection {Definitions and constructions}
In this section we define a 2-dimensional $\Cw$-complex corresponding to a symmetric SB-algebra $\Lambda$.
Associate with each $G$-cycle $z$ of length $k$ a $k$-gon $F_z$ with an oriented border. The sides of $F_z$ are labeled with the vertices of $Q_e$ which lie on $z$ (in the counter-clockwise order in the orientation of $F_z$). Consider a $\Cw$-complex $C=C(\Lambda)$ which is obtained from the resulting set of polygons by identifying oppositely oriented edges labeled by the same vertex. Since each vertex of $Q_e$ belongs to exactly two $G$-cycles, $C$ is an oriented manifold (without boundary).
\begin{definition}
$\Cw$-complex $C(\Lambda)$ is called {\it Brauer complex} of $\Lambda$.
\end{definition}
Denote by $\Gamma=\Gamma(\Lambda)$ the Brauer graph of $\Lambda$. For a vertex $V\in V(\Gamma)$, consider a cyclic permutation $\pi_V$ of half-edges, incident with $V$, which is defined by passing along the corresponding $A$-cycle. A 'picture' of a graph $\Gamma$ on an oriented surface also determines, for any vertex of the graph, a cyclic permutation on the set of incident half-edges, which agrees with orientation. There exists an embedding $i_\Gamma$ of $\Gamma$ into an oriented surface $M$, which preserves the cyclic permutations ($i_\Gamma$ and $M$ are uniquely defined up to a homeomorphism). Note that we consider {\it strict} embeddings, i.e. such embeddings that each connectivity component of $M \setminus \Gamma$ is homeomorphic to an open disk). See~\cite{7} for the construction of embedding. It follows from the construction of embedding that the connectiity components of $M \setminus \Gamma$ correspond to the $G$-cycles of $\Lambda$. Now it is clear that $M$ is a geometric realization of $C(\Lambda)$ and that the 1-skeleton $S_M$ of $C(\Lambda)$ is isomorphic as a graph to $\Gamma$ (we will refer to $S_M$ as $\Gamma$). In particular, the vertices (edges) of $C(\Lambda)$ are in one-to-one correspondence with the $A$-cycles (resp., vertices) of $Q_e$. It is to be mentioned that the arrows of $Q_e$ are in one-to-one correspondence with the angles of the 2-dimensional faces of $C(\Lambda)$.
\begin{definition}
{\it Perimeter} of a 2-dimensional face of $C(\Lambda)$ is the number of its edges, taking multiplicities into account (i.e., perimeter is the length of the corresponding $G$-cycles).
\end{definition}
\subsection{Invariants of stable equivalence}
Observe that $C(\Lambda)$ is an oriented surface. The following statement holds since the Euler characteristic of an oriented surface is even.
\begin{proposition}
If in the extended quiver $Q_e$ of $\Lambda$ the number of $A$-cycles is $k$, the number of $G$-cycles is $g$ and the number of vertices is $n$, then $k+g-n$ is even.
\end{proposition}
\begin{remark}
This statement was proved in~\cite{2} without topological arguments (Lemma 3.2).
\end{remark}
\begin{definition} 
The value $k+g-n$ is called the {\it genus} of $\Lambda$ (and of $C(\Lambda)$).
\end{definition}

\noindent In~\cite{2} it is proved that the multiset of lengths of $G$-cycles, as well as the number of $A$-cycles, is invariant under stable equivalence. By Rickard's Theorem, the derived equivalence of self-injective algebras implies stable equivalence (See~\cite{6}). The number of isomorphism classes of simple modules (i.e., the number of vertices of $Q_e$) is also stable invariant (See~\cite{4}). Therefore we get

\begin{proposition}
The multiset of perimeters of faces, the number of vertices and the genus of
$C(\Lambda)$ are invariant under derived equivalence.
\end{proposition}

\noindent It was shown in~\cite {2} that the free rank of the Grothendieck group of the stable category
$\text{stmod-}\Lambda$ equals $n-k$ if and only if $\Gamma(\Lambda)$ is not bipartite. Therefore, we have

\begin{proposition} \label{bipartite}
Derived (stable) equivalence preserves the property of the Brauer graph to be bipartite.
\end{proposition}

It should be mentioned that for algebras of genus $0$ this invariant gives nothing new, since an embedded into a sphere graph is bipartite if and only if the perimeters of all its faces are even. But there are algebras of genus 1, the derived categories of which are not distinguished by the previously discussed invariants, but which are not equivalent by Proposition~\ref{bipartite}.

\begin{example}
Consider the following symmetric SB-algebras $\Lambda_1$ and $\Lambda_2$:

\smallskip
\noindent 
\begin{tabular}{rp{12cm}}
$\xymatrix{
&&3 \ar@<2pt>[ddl]^\gamma \ar@<-1pt>[ddl]_\eta \\
&&&\\
&1 \ar@<2pt>[rr]^\alpha \ar@<-1pt>[rr]_\delta &&2 \ar@<1pt>[uul]^\varepsilon \ar@<-2pt>[uul]_\beta 
}$ & The quiver $Q_1$ of $\Lambda_1$ consists of vertices $1,2,3$ and arrows $$\alpha, \delta: 1\to 2, \text{ } \beta, \varepsilon: 2\to 3, \text{ } \gamma, \eta:3\to 1.$$ Ideal $I_1$ of relations of $\Lambda_1$ is generated by the elements $$\alpha\beta, \text{ }\beta\gamma,\text{ }\gamma\delta, \text{ }\delta\varepsilon,\text{ }\varepsilon\eta,\text{ }\eta\alpha,\text{ }
\alpha\epsilon\gamma-\delta\beta\eta,\text{ }\varepsilon\gamma\alpha-\beta\eta\delta,\text{ }
\gamma\alpha\varepsilon-\eta\delta\beta.$$
\end{tabular}

\noindent\begin{tabular}{p{12cm}r}
\noindent The quiver $Q_2$ of $\Lambda_2$ consists of vertices $1, 2, 3$ and arrows $$\alpha_1:1\to 2,  \text{ } \beta_1: 2\to 3, \text{ } \gamma_1:3\to 1, \text{ } \delta_1:1\to 3, \text{ } \varepsilon_1: 3\to 2, \text{ } \eta_1:2\to 1.$$ Ideal $I_2$ of relations of $\Lambda_2$ is generated by the elements
\begin{multline*}\alpha_1\beta_1, \text{ } \beta_1\gamma_1, \text{ } \gamma_1\delta_1, \text{ } \delta_1\varepsilon_1, \text{ } \varepsilon_1\eta_1, \text{ } \eta_1\alpha_1, \text{ }
\alpha_1\eta_1\delta_1\gamma_1-\delta_1\gamma_1\alpha_1\eta_1, \\\text{ } \beta_1\varepsilon_1-
\eta_1\delta_1\gamma_1\alpha_1, \text{ } \varepsilon_1\beta_1-\gamma_1\alpha_1\eta_1\delta_1.\end{multline*}
& $\xymatrix{
&&3 \ar@<2pt>[ddl]^{\gamma_1} \ar@<2pt>[ddr]^{\varepsilon_1} \\
&&&\\
&1 \ar@<2pt>[rr]^{\alpha_1} \ar@<1pt>[uur]^{\delta_1} &&2 \ar@<1pt>[ll]^{\eta_1} \ar@<1pt>[uul]^{\beta_1} 
}$
\end{tabular}

\noindent It is easy to see that $\Lambda_1$ and $\Lambda_2$ are algebras with 3 simple modules, with one $G$-cycle of length $6$
($(\alpha\beta\gamma\delta\varepsilon\eta)$ and  $(\alpha_1\beta_1\gamma_1\delta_1\varepsilon_1\eta_1)$, respectively) and with 2 $A$-cycles of multplicities $1$ ($c^1_1=(\alpha\varepsilon\gamma)$, $c^1_2=(\delta\beta\eta)$ and $c^2_1=(\varepsilon_1\beta_1)$, $c^2_2=(\gamma_1\alpha_1\eta_1\delta_1)$).
In particular, $\Lambda_1$ and $\Lambda_2$ have genus 1.
But $\Gamma(\Lambda_1)$ is bipartite (it consists of 2 vertices, connected by 3 edges) whereas $\Gamma(\Lambda_2)$ is not (the edge, corresponding the vertex 1 of $Q_2$ is a loop). Therefore, $\Lambda_1$ and $\Lambda_2$ are not derived equivalent.
\end{example}

Despite existence of an 'additional' invariant, the invariants and equivalences which are discussed in this paper are not enough to classify algebras of positive genus, in contrast to the 'spherical' case, which is treated in section 5 (see also example~\ref{example}).

\begin{proposition}\label{bijection} Correspondence $\Lambda\mapsto C(\Lambda)$ gives a bijection from the set of
(pairly non-isomorphic) indecomposable symmetric $\Sb$-algebras to the set of (pairly non-isomrphic) pairs $(C$, $f)$, where

\begin{enumerate}
	\item $C$ is a $\Cw$-complex homeomorphic to 2-dimensional oriented manifold with fixed orientaton;
  \item $f$ is an arbitrary map from the $0$-skeleton of $C$ to $\mathbb{N}$. 
\end{enumerate}
\end{proposition}

\begin{proof}
It remains to show that a Brauer complex uniquely determines a symmetric $\Sb$-algebra. It follows from the fact the
1-skeleton of Brauer complex has a structure of Brauer graph, which uniquely determines a symmetric SB-algebra (see~\cite{1})\footnote{ In~\cite{1} it was shown that a symmetric SB-algebra is uniquely determined by the (labeled) Brauer graph and certain parameters. It can be easily shown that these parameters are excessive and can be eliminated.}.
\end{proof}

\section{Elementary tilting complexes}
\subsection{Definition of elementary tilting complex}
Fix an edge $i$ of $C$ (equivalently, fix a vertex $i$ in quiver $Q_e$), and suppose that there are other edges in $C$.
We distinguish three cases.
\begin{enumerate}
    \item $i$ is a leaf of $\Gamma$. Equivalently, in the quiver $Q_e$ there is a loop $\alpha_i$ at vertex $i$ and this loop is an $A$-cycle (i.e., it annihilates all other arrows of $Q_e$).
    \item $i$ is a loop, which bounds some face of $C$. Equivalently, in the quiver $Q_e$ there is a loop $\alpha_i$ at vertex $i$ and this loop is a $G$-cycle. In this case there is a unique $A$-cycle passing through $i$ (this cycle contains at least 3 arrows, one of which is $\alpha_i$).
    \item For $r=1,2$ the end $C_{i,r}$ of the edge $i$ is incident with an edge $i_r \neq i$, such that $\pi_{C_{i,r}}(i_r)=i$. We permit $i_1=i_2$ and we permit $i$ to be a loop (i.e., $C_{i,1}=C_{i,2}$). Equivalently, there is no loop at vertex $i$ of $Q_e$, i.e. the vertices $i_1, i_2$ which precede $i$ on both $A$-cycles passing through $A$ ($c_{i,1}$ and $c_{i,2}$) are different from $i$. 
\end{enumerate}
In each of these cases, to the edge $i$ we put in correspondence a complex $T_i$ as follows. For a vertex $j\in V(Q_e)$ we denote by $P_j$ the indecomposable left projective $\Lambda$-module, which corresponds to $j$. For $i \neq j$, denote by $T_{ij}$ the complex $\dots\to 0\to P_j\to 0\to\dots$ concentrated in degree $0$. If $i$ is a leaf of $\Gamma$, define complex $T_{ii}$ by
$$
T_{ii}:\dots\to 0\to P_j\stackrel{\beta_i}{\longrightarrow} P_i \to 0\to \dots
$$
where $j \in V(Q_e)$, $j\neq i$ is the vertex preceding vertex $i$ on the (unique) $G$-cycle, which contains $i$; $\beta_i \neq \alpha_i$ is the arrow preceding $\alpha_i$ on the same $G$-cycle.

\noindent If $i$ is a loop which bounds some face of $C$, define $T_{ii}$ by
$$
T_{ii}:\dots\to 0\to P_j\bigoplus P_j\xrightarrow{(\beta_i,\text{ }\beta_i\alpha_i)}
P_i \to 0\to \dots
$$
where $j \in V(Q_e)$, $j\neq i$ is the vertex preceding vertex $i$ on the (unique) $A$-cycle, which contains $i$; $\beta_i \neq \alpha_i$ is the arrow preceding $\alpha_i$ on the same $A$-cycle.

\noindent Otherwise, define $T_{ii}$ by
$$
T_{ii}:\dots\to 0\to P_{i_1}\bigoplus P_{i_2}\xrightarrow
{(\beta_{i}^1,\text{ }\beta_{i}^2)} P_i \to 0\to \dots
$$
where $i_1, i_2$ are the vertices preceding $i$ on the $A$-cycles $c_{i,1}$ and $c_{i,2}$, respectively; $\beta_i^1,\beta_i^2$ are the respective arrows preceding $\alpha_i$.
Finally, set $T_i=\bigoplus_{j=1}^{n}T_{ij}$.

\begin{proposition}
$T_i$ is a tilting complex over $\Lambda$.
\end{proposition}
\begin{proof}
We verify that $T_i$ satisfies the two conditions from the definition of tilting complex. In the definition of $T_i$ we distinguished three cases. We show verification only for the third case, the other cases are treated in the same way.

First, we must verify that $D^b(\Lambda)= Add(T_i)$, where $Add(T_i)$ is the smallest triangulated subcategory, which contains all direct summands of object $T_i$. It is enough to verify that all objects of the form $0\to P_j\to 0$ belong to $Add(T_i)$. For $i\neq j$ this is by definition of $T_i$. For $i=j$ it is easy to see that $P_i[-1]$ is the third term of the triangle, which corresponds to the natural embedding of  $T_{ii_1}\bigoplus T_{ii_2}$ into $T_{ii}$. It follows that $T_i$ satisfies the first condition.

Now we verify that $\text{Hom}_{D^b(\Lambda)}(T_i,T_i[r])=0$ for $r\in \mathbb{Z}\setminus0$. It is enough to proof that for each $j \in V(Q_e)$ $\text{Hom}_{D^b(\Lambda)}(T_{ii},T_{ij}[-1])=\text{Hom}_{D^b(\Lambda)}(T_{ij}[-1],T_{ii})=0$. Each morphism from $T_{ij}$ to $T_{ii}$ is determined by a morphism $f:P_j\to P_i$, where $f$ is a multiplication by a linear combination of paths with starting point $j$ and endpoint $i$. Each of these paths ends either with $\beta_{i}^1$ or with $\beta_{i}^2$. Therefore $f$ factors through $(\beta_{i}^1,\beta_{i}^2):P_{i_1}\bigoplus P_{i_2}\to P_i$. It follows that $f$ is homotopic to zero. Similarly, each morphism from $T_{ii}$ to $T_{ij}$ is determined by a morphism $f:P_i\to P_j$, where $f$ is a multiplication by a linear combination $S$ of paths with starting point $i$ and endpoint $j$. Suppose that $S$ has nonzero summands. Since $i\neq j$, the underlying paths are not maximal. Multiplying $S$ by $\beta_{i}^1$ or by $\beta_{i}^2$ from the left, we again get a nontrivial sum of linearly independent summands. This contradicts the definition of morphism of complexes. Therefore $f=0$ and $\text{Hom}_{D^b(\Lambda)}(T_{ii},T_{ij}[-1])=0$.
\end{proof}

\subsection{Elementary transformations of Brauer complexes}
Now we define {\it elementary transformations} of Brauer complexes. We will prove below that in terms of algebras, an elementary transformation puts an algebra $\Lambda$ to the endomorphism algebra of one of the above defined tilting complexes over $\Lambda$. We fix convention that under elementary transformation the vertices are fixed, the configuration of edges (labeled with vertices of a quiver) --- and therefore the configuration of  faces (labeled with $G$-cycles) --- is changed. In other words, we identify the edges (and faces) by their labels, not by the vertices incident to them. The pictures below illustrate the simplest cases, in general they can be quite different.

\begin{definition}
Let $C$ be a Brauer complex, let $Q_e$ be the corresponding extended quiver. Let $a \in E(C)$, $V\in V(C)$, let $F$ be a face of $C$. Permutations $\Next_F$: $V(C) \rightarrow V(C)$ and $E(C) \rightarrow E(C)$ are induced by the counter-clockwise order of vertices and edges in the orientation of $F$. Recall that $\pi_V$ denotes the permutation of half-edges incident with vertex $V\in V(C)$, which is defined by passing along the corresponding $A$-cycle $v$. By abuse of language, we will name half-edges after correspondent edges. Thus by abuse of language for a loop $a$ both situations $\pi_V(a)=a$ and $\pi_V(a)\neq a$ can happen. However, from the context it will always be clear which half-edge is meant. 
\end{definition}

\subsubsection{Transformation of type 1: shift of a leaf}
\begin{figure}[ht]
\begin{center}
		\includegraphics{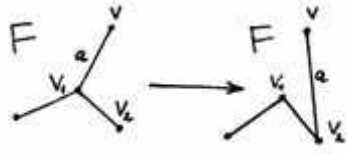} 	
		\caption{Shift of a leaf}
  	\label{fig:leaf}
	\end{center}
\end{figure}
Let $V \in V(C)$ be a dangling vertex. Suppose that the edge (the face) incident with $V$ is labeled by $a$ (resp., by $F$). 
Let $V_1$ be the second vertex incident with $a$. Put $V_2 = \Next_F(V_1)$, $a_1= \Next_F(a)$. Now shift edge $a$, so that $a$ becomes incident with $V$ and $V_2$ and $a=\Next_F(a_1)$. 

\subsubsection{Transformation of type 2: shift of a loop}
Let $a$ be a loop at vertex $V_1$, bounding some face $F_1$. Let $F_2$ be the second face, incident with $a$, put $V_2 = Next_{F_2}(V_1)$, $a_1= Next_{F_2}(a)$. Replace loop $a$ with a loop at vertex $V_2$, which lies inside $F_2$ after $a_1$. Note that $F_1$ is again bounded by a loop, which separates it from $F_2$.
\begin{figure}[ht]
\begin{center}
		\includegraphics{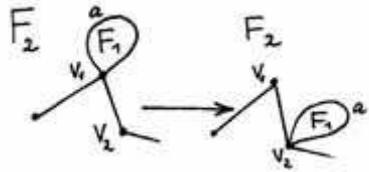}
  	\caption{Shift of a loop}
  	\label{fig:loop}
\end{center}
\end{figure}
\subsubsection{Transformation of type 3: the general case} \label{mainshift}
Let $a$ be an edge. Suppose that the vertices (faces) incident with $a$ are labeled by $V_1, V_2$ (resp., by $F_1$ and $F_2$; we permit $F_1=F_2$). 
For $i=1,2$ put ${V_i}'=Next^{-1}_{F_i}(V_i)$, $a_i=Next^{-1}_{F_i}(a)$. Shift $a$ so that it becomes incident with ${V_1}'$ and ${V_2}'$, separates $F_1$ from $F_2$ and lies after $a_i$ on the new boundary of $F_{3-i}$.
\begin{figure}[ht]
\begin{center}
		\includegraphics{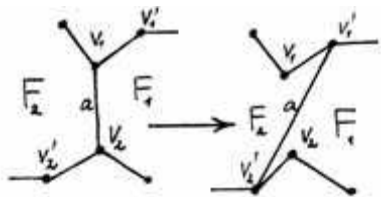}
  	\caption{The general case}
  	\label{fig:general}
\end{center}
\end{figure}
\begin{definition}
We call the transformations of types 1-3 {\it tilting transformations}. The resulting complex is denoted by $C(a)$.
\end{definition}

\subsection{Correspondence}
\begin{proposition} \label{correspondence}
Let $\Lambda$ be an $\Sb$-algebra, $C=C(\Lambda)$, $a\in E(C)$. Let $T_a$ be the tilting complex which corresponds to $a$. Then $End_{D^b(\Lambda)} T_a$ is a symmetric $\Sb$-algebra with Brauer complex $C(a)$ ($C$ and $C(a)$ have the same multiplicities of vertices).
\end{proposition}
\begin{proof} 
Denote by $Q_e$ the extended quiver of $\Lambda$. By Rickard's theorem \cite{5}, $\Lambda_a=End_{D^b(\Lambda)} T_a$ is derived equivalent to $\Lambda$. Since $\Lambda$ is a symmetric algebra, $\Lambda_a$ is a symmetric algebra, too. By Pogorjaly's result, an algebra, which is stable equivalent to an SB-algebra, is an SB-algebra, too \cite{4}. Therefore, by another Rickard's theorem \cite{6} $\Lambda_a$ is an $SB$-algebra. Let $e=\sum_1^n e_i$ be the decomposition of unity of $\Lambda_a$, which corresponds to the decomposition
$T_a=\bigoplus _{i=1}^n T_{ai}$. Since the number of simple modules is invariant under derived equivalence, $\Lambda_a$ is an algebra with $n$ simple modules and therefore $\{e_i\}$ is a set of primitive orthogonal idempotents. Set $f_a=1-e_a\in\Lambda$, and denote $\Lambda_{-a}=f_a\Lambda f_a$. Since for $i\neq a$ the complexes $T_{ai}$ are concentrated in degree $0$, we have $\Lambda_{-a}=\text{End}_{\Lambda}\bigoplus_{i\neq a}P_i=\text{End}_{D^b(\Lambda)}
\bigoplus_{i\neq a}T_{ai}$. Consider Brauer complex $C_{-a}$, obtained from $C$ by deletion of an edge $a$ (if $a$ is a leaf, we delete it with the incident dangling vertex). The marks on the remaining vertices are preserved.
We need the following lemma.

\begin{lemma} \label{lemma}
The symmetric $\Sb$-algebra which corresponds to $C_{-a}$ is isomorphic to
$\Lambda_{-a}$.
\end{lemma}
\begin{proof}
We consider the case when $C(a)$ is obtained from $C$ by a transformation of type 3 (i.e., $a$ is not a loop which bounds a face and not a leaf). The other cases are treated in the same way. Denote the arrows of $Q$ incident with $a$ by $\alpha,\beta,\gamma,\delta$, so that $\alpha\beta\neq 0$ and $\gamma\delta\neq 0$. The elements of $\Lambda_{-a}$ are linear combinations of paths whose starting points and endpoints differ from $a$. It is clear that $\Lambda_{-a}$ is generated as algebra by idempotents $e_i$, where $i\neq a$, by arrows of $Q$ different from $\alpha,\beta,\gamma,\delta$ and by the elements $\alpha\beta,\gamma\delta$. Observe that in terms of quivers $\Lambda_{-a}$ can be obtained from $\Lambda$ in the following way: the arrows $\alpha$ and $\beta$, lying on a common $A$-cycle, are replaced with an arrow $\alpha\beta$ on the same $A$-cycle (respectively, the arrows $\gamma$ and $\delta$ are replacesd with an arrow $\gamma\delta$). This implies the claim.
\end{proof}
\noindent We return to the proof of proposition~\ref{correspondence}. Observe that the symmetric $\Sb$-algebra which corresponds to $C(a)_{-a} = C_{-a}$ is isomorphic to $\Lambda_{-a}$. To obtain the Brauer complex of $\Lambda_a$ from $C_{-a}$ we need to add an edge on some face of $C_{-a}$ (the multiplicities of vertices are preserved). It should be noted that all arrows of the quiver of $\Lambda_{-a}$ except at most two coincide with the respective arrows of the quiver of $\Lambda_{a}$. The arrows which don't coincide, are products of two or three arrows of the quiver of $\Lambda_{a}$. Again, we finish the proof only for the case when $C(a)$ is obtained from $C$ by tilting transformation of type 3; the other cases are treated in the same way. 
For $i=1,2$ denote by $b_i$ the edge, which precedes $a_i$ on $F_i$ in counter-clockwise order, i.e. $b_i$ precedes $a_i$ on a $G$-cycle (see notations in~\ref{mainshift}).
Denote by $\mu$ (by $\rho$) the arrow in $Q_e$ which corresponds to the angle at vertex $V_1'$ included between $a_1$ and $b_1$ (resp., to the angle at $V_2'$ included between $a_2$ and $b_2$).
Define elements
$\alpha_1,\beta_1,\gamma_1,\delta_1\in End_{D^b(\Lambda)} T_a$
such that
$\alpha_1\beta_1=\mu$, $\gamma_1\delta_1=\rho$
in
$\text{End}_{D^b(\Lambda)}\bigoplus_{i\neq a}T_{ai}$.
Each of these elements is induced by a morphism between two indecomposable summands of $T_a$:

\noindent$
\alpha_1:
$
$$
\begin{CD}
\dots @>>> 0 @>>> P_{b_1} @>>> 0 @>>> \dots @.\\
 @. @VVV @VV\binom{\mu}0 V @VVV @. @.\\
\dots @>>> 0 @>>> P_{a_1}\bigoplus P_{a_2} @>(\alpha,\gamma)>> P_a @>>> 0 @>>> \dots
\end{CD}
$$
$
\beta_1:
$
$$
\begin{CD}
\dots @>>> 0 @>>> P_{a_1}\bigoplus P_{a_2} @>(\alpha,\gamma)>> P_a @>>> 0 @>>> \dots\\
@. @VVV @VV(id,0) V @VVV @. @.\\
\dots @>>> 0 @>>> P_{a_1} @>>> 0 @>>> \dots @.\\
\end{CD}
$$
$
\gamma_1:
$
$$
\begin{CD}
\dots @>>> 0 @>>> P_{b_2} @>>> 0 @>>> \dots @.\\
 @. @VVV @VV\binom0{\rho} V @VVV @. @.\\
\dots @>>> 0 @>>> P_{a_1}\bigoplus P_{a_2} @>(\alpha,\gamma)>> P_a @>>> 0 @>>> \dots
\end{CD}
$$
$
\delta_1:
$
$$
\begin{CD}
\dots @>>> 0 @>>> P_{a_1}\bigoplus P_{a_2} @>(\alpha,\gamma)>> P_a @>>> 0 @>>> \dots\\
@. @VVV @VV(0,id) V @VVV @. @.\\
\dots @>>> 0 @>>> P_{a_2} @>>> 0 @>>> \dots @.\\
\end{CD}
$$
The elements $\alpha_1,\beta_1,\gamma_1,\delta_1$ are not invertible, since for $i\neq a$ $H^*(T_{aa})\neq H^*(T_{ai})$. Therefore these are the arrows $\mu$ and $\rho$ (in $\Lambda_{-a}$) which are products of two arrows of $\Lambda_{a}$. Now observe that in terms of Brauer complexes, transformation of the quiver of $\Lambda_{-a}$ to $\Lambda_{a}$ is insertion of edge labeled by $a$, incident with $V_1'$ and $V_2'$ into the union of faces $F_1$ and $F_2$.
\end{proof}
\begin{corollary}
Let $\Lambda_1$ and $\Lambda_2$ be symmetric $\Sb$-algebras, let $C_1$ and $C_2$ be their Brauer complexes. Suppose that $C_2$ can be obtained from $C_1$ by a sequence of tilting transformations. Then $\Lambda_1$ and $\Lambda_2$ are derived equivalent.
\end{corollary}
\begin{proof}
The statement follows from Proposition~\ref{correspondence}, Lemma~\ref{lemma} and the Rickard's Theorem.
\end{proof}
\begin{example} \label {example}
Consider decagons $D_1$ and $D_2$. Fix an orientation on each of decagons. Mark the edges of $D_1$ (of $D_2$) with letters $a$, $b$, $c$, $d$, $e$ so that they form a word $abcdeabcde$ (resp., $abcdeadebc$) in counter-clockwise order. In each decagon, identify the edges which are marked by the same letter in such way that the resulting manifolds are oriented. It's easy to see that both complexes (we call them $C_1$ and $C_2$) have 2 vertices, 5 edges, one face, i.e. they are homeomorphic to a sphere with two handles. Moreover, the 1-skeletons of $C_1$ and $C_2$ are bipartite graphs. But these complexes cannot be obtained from each other by tilting transformations: any complex $C'$, obtained from the complex $C_1$, is isomorphic to $C_1$. This construction gives pairs of symmetric $\Sb$-algebras of genus 2, for which the methods given in present paper are not enough to determine whether they are derived equivalent or not.
\end{example}
\section{Algebras of genus 0}
Now we prove that if Brauer complex of $\Lambda$ is homeomorphic to a sphere, then the multiset of perimeters of its faces and the multiset of multiplicities of vertices determine the class of derived equivalence of $\Lambda$. For a start, we don't take into consideration the multiplicities of vertices, i.e. we consider graphs with non-labeled vertices. We fix plane graphs $\Gamma_1$ and $\Gamma_2$ with the same multisets of perimeters of faces and show that $\Gamma_2$ can be obtained from $\Gamma_1$ by a sequence of tilting transformations (statements from Lemma~\ref{lemma1} to Proposition~\ref{itog1}).
\begin{definition}
Graphs which can be obtained from each other by a sequence of tilting transformations will be called {\it chain equivalent} graphs.
\end{definition}
\begin{lemma} {\label {lemma1}}
Let $\Gamma$ be a plane graph, $A \in V(\Gamma)$. There exists a plane graph $\Gamma'$, chain equivalent to $\Gamma$, in which the vertex $A$ is incident with all edges and one of the following conditions holds:
\begin{enumerate}
    \item $\Gamma'$ has no loops
    \item Each edge of $\Gamma'$ is either a leaf or a loop at vertex $A$ (i.e., there are no multiedges in $\Gamma'$ except for loops).
\end{enumerate}
\end{lemma}
\begin{definition}
Plane graph of this form is called a {\it reduced graph}. 
\end{definition}
\begin{proof}
Consider among graphs, which are chain equivalent to $\Gamma$, a graph $\Gamma'$ with a maximal degree of $A$. Observe that all edges of $\Gamma'$ are incident with $A$. Indeed, otherwise there are vertices $B, C \neq A$ and an edge $e\in E(B,C)$ such that either $B$ or $C$ is incident with $A$ (without loss of generality, $B$) and such that the edge $\pi_B(e) \in E(A,B)$. If $B\neq C$, we apply to $e$ a transformation of type 3. If $B=C$, we apply to $e$ a transformation of type 2 so that $e$ shifts from $B$ to $A$. Thus the degree of $A$ can be increased, a contradiction. It follows that there are three types of edges in $\Gamma'$:
\begin{enumerate}
    \item[a)] a loop at vertex $A$;
    \item[b)] edges which form a multiedge incident with $A$; 
    \item[c)] a leaf $(A,X)$.
\end{enumerate}
For further convenience, elements of type a) don't belong to type b). We show that in $\Gamma'$ edges of types a) and b) cannot exist simultaneously. Suppose that there is a loop $a$, leaves $a_1=\pi_A(a)$, $a_2=\pi_A(a_1), \dots$, $a_s=\pi_A(a_{s-1})$ and an edge $b=\pi_A(a_s)$ of type b). Consider the edge $c=\pi_B(b)$. By transformations of type 1, we shift $a_1, \dots, a_s$ along $a$. Now there are no edges between $a$ and $b$ around $A$, and we can apply a transformation of type 3 to the edge $b$ and $b$ becomes a loop. This increases the degree of $A$, a contradiction.
\end{proof}
\begin{definition}
A reduced graph which has no loops is called a {\it reduced graph of type 1}. 
\end{definition}
\noindent Observe that the border of any face of a reduced graph of type 1 is formed by several pairs of edges $(A,B_1), \dots, (A,B_k)$ and by several leaves (any leaf is counted in the perimeter of the face twice). Observe that a reduced graph of type 1 is bipartite. 
\begin{definition}
A reduced graph which has loops is called a {\it reduced graph of type 2}. 
\end{definition}
\noindent In a reduced graph of type 2, any edge which is not a leaf is a loop. Observe that a reduced graph of type 2 is not bipartite. 

Let $\Gamma_1'$ and $\Gamma_2'$ be reduced graphs, chain equivalent to $\Gamma_1$ and to $\Gamma_2$, respectively. By Proposition~\ref{bipartite}, $\Gamma_1'$ and $\Gamma_2'$ are of the same type. We will show that all reduced graphs of the same type, with the same multisets of perimeters of faces, are chain equivalent:

\smallskip
\noindent {\bfseries I. Reduced graphs of type 1.} Fix a graph $\Gamma$ of type 1. 
\begin{lemma} \label{reduced1}
Each reduced graph $\Gamma$ of type 1 is chain equivalent to a reduced graph $\Gamma'$ of type 1, which has at most two non-dangling vertices.
\end{lemma}
\begin{proof}
Consider among reduced graphs, which are chain equivalent to $\Gamma$, a graph $\Gamma'$ with maximal number of dangling vertices. Let $A$ be the vertex of $\Gamma'$, which is incident with all edges. We show that $\Gamma'$ has at most two non-dangling vertices (including $A$). Indeed, let $b\in E(A,B)$ and $c\in E(A,C)$ be two edges of type 2 ($B\neq C$) such that there are only leaves between $b$ and $c$ in clockwise order around A. As above, by transformations of type 1 we obtain a graph, in which there are no leaves between $b$ and $c$ (around $A$). Suppose that $C$ has degree $2$. Applying the transformation of type 3 to $c$ (shift along $b$), we get a reduced graph with a greater number of leaves, since $C$ becomes a leaf. In order to transform $C$ to a leaf when deg$(C)=r$, we need to carry out the same operations with $r-1$ edges, which are incident with $C$.
\end{proof}
Consider a reduced graph $\Gamma'$ which was obtained in lemma~\ref{reduced1}.
It is easy to see that the faces of $\Gamma'$ and the edges of $\Gamma'$ which are not leaves can be cyclically numbered by $1, 2...\dots, g$ so that the border of the face number $i$ consists of the edges number $i$ and $i+1$ and several inner leaves. It should be mentioned that if $g=1$ then $\Gamma'$ is a tree in a form of star, and we get Brauer trees, which were studied by Rickard in~\cite{6}, as a first application of the criterion of derived equivalence.

\begin{lemma}\label{canonical}
Graph $\Gamma'$ is chain equivalent to a graph of the same form (i.e., as in lemma~\ref{reduced1}), in which the perimeters of faces are in ascending ordering.
\end{lemma}
\begin{proof}
It's enough to show how to 'transpose' two faces, see Figure~\ref{fig:circles}.
\begin{figure}[ht]
\begin{center}
		\includegraphics{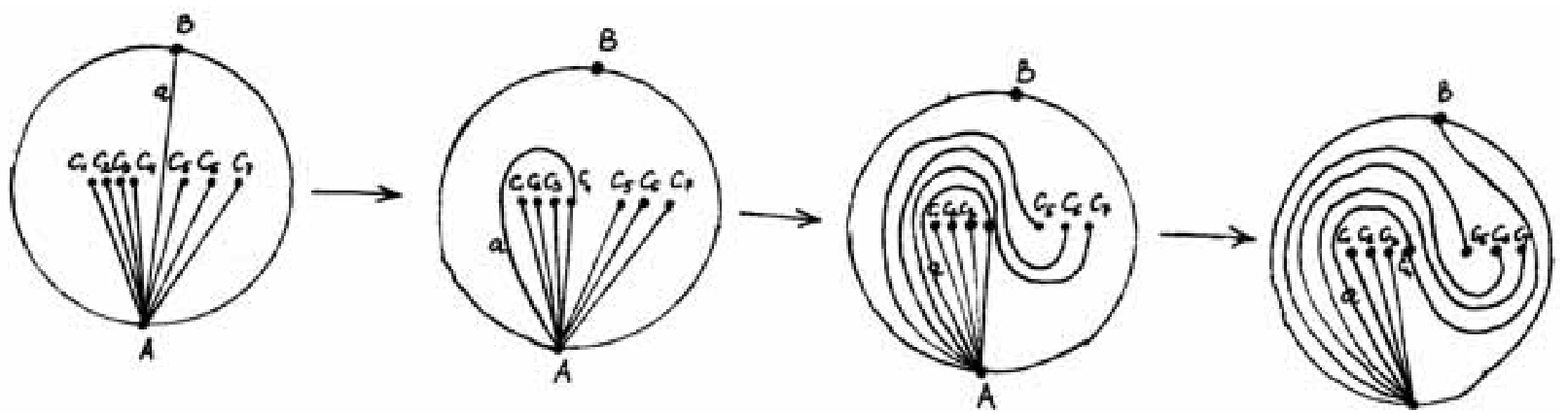}
		\caption{to Lemma~\ref{canonical}}
  	\label{fig:circles}
\end{center}
\end{figure}
\end{proof}
We see that any bipartite plane graph is chain equivalent to a (unique) {\it canonical representative} (we will also say "a {\it graph in canonical form}") --- a graph in which the perimeters of faces are in ascending ordering. Two graphs with the same multisets of perimeters are chain equivalent to the same canonical representative, and therefore they are chain equivalent to each other.

\smallskip
\noindent {\bfseries II. Reduced graphs of type 2.}

Consider a reduced graph $\Gamma$ of type 2.
First suppose that $A$ is the only vertex of $\Gamma$, i.e. {\bf all edges of $\Gamma$ are loops} and $n=g-1$, where $n$ is the number of vertices of $Q_e$ and $g$ is the number of $G$-cycles. 
Consider a graph $T=T(\Gamma)$, which is plane dual to $\Gamma$. $T$ is a tree with $g-1$ edges and $g$ vertices.
Observe that the transformations of type 1 cannot be applied to $\Gamma$. The transformations of types 2 and 3 can be described in terms of $T$ as follows.

\begin{figure}[ht]
\begin{center}
		\includegraphics{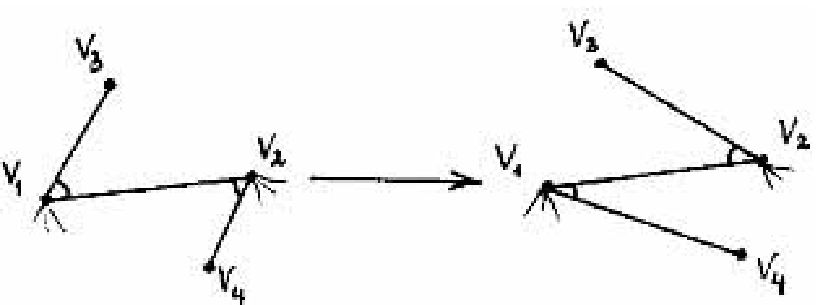}
		\caption{Flip}
  	\label{fig:flip}
\end{center}
\end{figure}

\begin{itemize}
	\item Transformation of type 2. A leaf $V_1V_2$ of $T$ (with dangling vertex $V_1$) is shifted around $V_2$ in arbitrary way. This transformation of a plane labeled tree will be called a {\it flip-over}.
	\item  Transformation of type 3. Suppose that $\pi^{-1}_{V_1}(V_1V_2)=V_1V_3$ and that $\pi^{-1}_{V_2}(V_1V_2)=V_2V_4$. Replace edges $V_1V_3$ and $V_2V_4$ with edges $V_1V_4$ and $V_2V_3$ in a way that $\pi_{V_1}(V_1V_2)=V_1V_4$ and $\pi_{V_2}(V_1V_2)=V_2V_3$. This transformation of a plane labeled tree will be called a {\it flip} (see Figure~\ref{fig:flip}; an arc between two edges in the pictures denotes absence of other edges). 
\end{itemize}

\begin{definition}
Plane trees with labeled vertices, which can be obtained from each other by flips and flip-overs, are called {\it equivalent}. Clearly, equivalent trees are dual to chain-equivalent graphs. 
\end{definition}
\begin{proposition} \label{equivalent}
Two plane trees with the same multisets of labeled vertices and the same degrees of correspondent vertices are equivalent.
\end{proposition}
\noindent We need the following lemma. 
\begin{lemma} \label{incident}
Let $V_1V_2$ be a leaf in a plane tree $T$ with dangling vertex $V_1$. 
Let $V_1,V_2,\dots, V_r$ be a path in $T$ such that $V_r$ is an non-dangling vertex. Then $T$ is equivalent to a tree, in which $V_1$ is adjacent with $V_r$. 
\end{lemma}
\begin{proof}
The proof is by induction on $r$. For $r=2$ the claim is trivial. Suppose that there is a number $i\in\{1,\dots,r\}$ such that deg$(V_i)\geq 3$. Consider the minimal such $i$. Without loss of generality we assume that $V_{i+1}\neq V$, where $V$ is such vertex that $\pi_{V_i}(V_iV_{i-1})=V_iV$.
If $i\neq 2$, replace edges $V_{i-1}V_{i-2}$ and $V_iV$ with $V_{i}V_{i-2}$ and $V_{i-1}V$ by a flip.
Otherwise, we make $V_2V_3$ follow $V_2V_1$ by several flip-overs, and then make the above flip. 
The distance between $V_1$ and $V_r$ decreases, and we apply the inductive hypothesis. 
If $i$ cannot be defined, consider the unique vertex $V_{r+1}\neq V_{r-1}$ adjacent with $V_r$.
Replace $V_{r-1}V_{r-2}$ and $V_rV_{r+1}$ with $V_{r}V_{r-2}$ and $V_{r-1}V_{r+1}$ by a flip. 
Again, the distance between $V_1$ and $V_r$ is decreased, and we apply the inductive hypothesis. 
\end{proof}
Now we prove Proposition~\ref{equivalent}.
\begin{proof}
The proof is by induction on the number of vertices. For $g=1$ the claim is trivial. Let $T_1$ and $T_2$ be two plane trees with $g$ vertices. 
Let a dangling vertex $V$ be adjacent with $V_1$ in $T_1$ and with $V_2$ in $T_2$. 
By Lemma~\ref{incident}, we can replace $T_1$ with an equivalent tree $T_3$ in which $V$ is adjacent with $V_2$. 
Let $T_3^1$ and $T_2^1$ be the trees, obtained from $T_3$ and $T_2$ by removing $V$ with the corresponding edge. 
They have the same degrees of correspondent vertices, and therefore they are equivalent by inductive hypothesis. 
It remains to show that it is still possible to carry out the sequence of transformations, which puts $T_3^1$ to $T_2^1$, when edge $V_2V$ is not deleted. After these transformation we will be able to flip-over the edge $V_2V$ to the required place. 

Start to apply the above sequence of transformations to $T_3$. We can encounter difficulties in the following cases:
\begin{itemize}
	\item When in $T_3$ the edge $V_2V$ is between two subsequent edges (around $V_2$) of $T_1^3$ and doesn't allow to make a flip. We cope with this by an arbitrary flip-over of $V_2V$. 
	\item 
	If $V_2$ is a dangling vertex in $T_3^1$, incident with an edge $V_2V_3$, and in $T_3^1$ it is possible to make a flip-over of $V_2V_3$. In $T_3$ instead of this flip-over we make the following sequence of transformations (Figure~\ref{fig:krakozabra}). 
\end{itemize}
\begin{figure}[ht]
\begin{center}
		\includegraphics{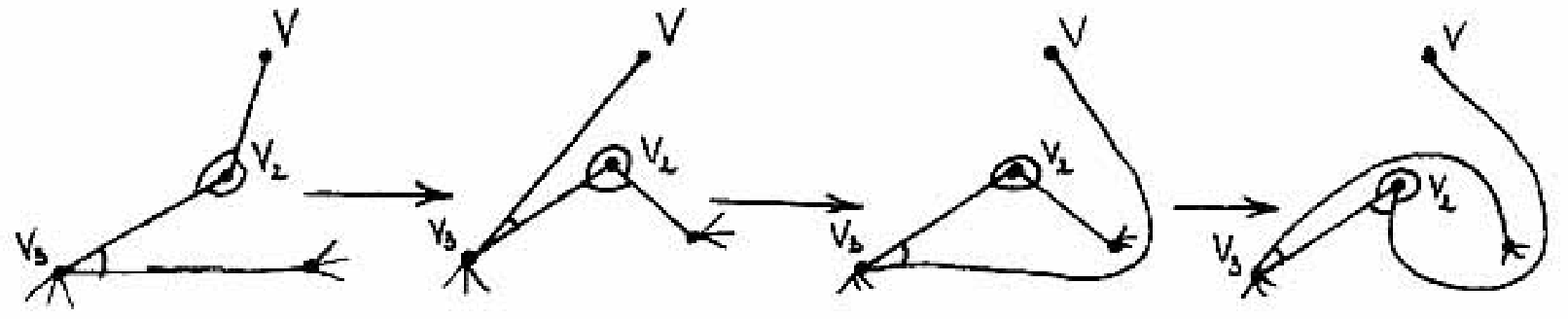}
  	\caption{to Proposition~\ref{equivalent}}
  	\label{fig:krakozabra}
\end{center}
\end{figure}	

\noindent This finishes the proof. 
\end{proof}
\noindent Now suppose that {\bf there are dangling vertices in $\Gamma$}. 
\begin{definition}
{\it External perimeter} of a face is the number of its edges, which separate it from other faces (in our case, these are loops). 
\end{definition} 
\begin{definition}
{\it Reduction} of graph $\Gamma$ is a graph $\mathcal{R}(\Gamma)$ which is obtained from $\Gamma$ by removing all dangling vertices. 
\end{definition}
\begin{proposition} \label{analogue}
Let $\Gamma_1$ and $\Gamma_2$ be reduced graphs of type 2. Suppose that there is a tilting transformation $p$ which puts $\mathcal{R}(\Gamma_1)$ to $\mathcal{R}(\Gamma_2)$. Suppose also that the correspondent labeled faces of $\Gamma_1$ and of $\Gamma_2$ have the same number of edges. 
Then graphs $\Gamma_1$ and $\Gamma_2$ are chain equivalent.
\end{proposition}
\begin{proof}
Let $l$ be the loop, which is shifted by $p$ and let $F_1$ and $F_2$ be the faces separated by $l$.
We need to obtain a sequence of transformations which would serve as an analogue of $p$ for $\Gamma_1$. Figure~\ref{fig:butterfly} illustrates the case when $F_1$ has inner leaves and $l$ is the only loop on the border of $F_1$. 
The case when there are other loops on the border of $F_1$ is even easier:  
the analogue of $p$ is a transformation of type 3, made after necessary flip-overs of leaves. 
\begin{figure}[ht]
\begin{center}
		\includegraphics{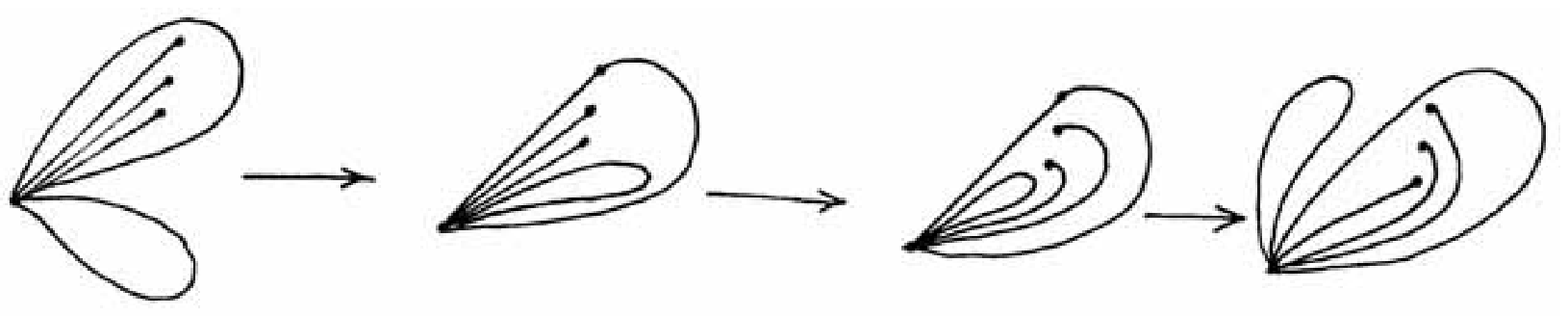}
  	\caption{to Proposition ~\ref{analogue}}
  	\label{fig:butterfly}
\end{center}
\end{figure}
\end{proof}

\begin{remark}\label{mainremark}
It follows from Propositions~\ref{equivalent} and~\ref{analogue} that the class of chain equivalence of a reduced graph of type 2 is determined by the multiset of pairs $(P(F_i)$, $p(F_i))$, where $P(F_i)$ is the perimeter and $p(F_i)$ is the external perimeter of the face $F_i$. 
\end{remark}
\begin{definition}
The multiset of pairs $(P(F_i)$, $p(F_i))$ will be called a {\it multiset of double perimeters} of graph $\Gamma$.
\end{definition}

\begin{proposition} \label{double}
Let $\Gamma_1$ and $\Gamma_2$ be reduced graphs of type 2 with the same multisets of perimeters of faces.
Then there exists a reduced graph $\Gamma_3$ of type 2, chain equivalent to $\Gamma_1$, such that the multisets of double perimeters of $\Gamma_2$ and $\Gamma_3$ are the same.
\end{proposition}

\begin{proof}
Let $\{(P_i$, $p_i)\}$ be the multiset of double perimeters of $\Gamma_2$, let $\{(P_i$, $p^1_i)\}$ be the multiset of double perimeters of $\Gamma_1$, for $i=1, \dots, g$.
Observe that $p_i\equiv P_i\equiv p^1_i\pmod 2$ for each $i\in \{1,\dots g\}$ and that $\sum_i p_i=2g-2=\sum_i p^1_i$.
Set $q_i=p_i^1$ for each $i$. 
Consider the following algorithm of 'transformation' of the multiset $\{q_i\}$ to the multiset $\{p_i\}$. Below we will show that for each step of this algorithm there is a chain equivalence of graphs, which properly changes their external perimeters. 

Consider maximal $k$ such that $q_i=p_i$ for all $i<k$.
\begin{enumerate}
	\item If $q_k<p_k$ then $q_j>p_j$ for some $j>k$. 
Replace $q_k$ with $q_k+2$ and replace $q_j$ with $q_j-2$. 
  \item 
Otherwise $q_k>p_k\geq 1$ and $q_j<p_j$ for some $j>k$. In this case we replace $q_k$ with $q_k-2$ and replace $q_j$ with $q_j+2$.
\end{enumerate}
Observe that at each step the number which is decreased is greater than two, so the resulting numbers are positive. 
Moreover, since $q_i \leq$ max$(p_i,p_i^1)$, at each step $q_i \leq P_i$ for all $i$.
Clearly, the multiset of numbers $q_i$ can be transformed to the multiset of numbers $p_i$ by these operations. 
To find the chain equivalences which correspond to these operations, we need the following lemma. 
\begin{lemma} \label{triv}
Let $T$ be a tree, let $V_1, V_2 \in V(T)$. If $V_1$ and $V_2$ are not both dangling vertices, then 
there exists a tree in which the degrees of all vertices are the same and the vertices $V_1$ and $V_2$ are adjacent.
\end{lemma}
\begin{proof}
The proof is by induction on the number of vertices in $T$. 
\end{proof}
\noindent We return to the proof of Proposition~\ref{double}.
We need a sequence of tilting transformations under which the multiset of external perimeters changes in accordance to the above algorithm. Suppose that we are to change the external perimeters $q_i$ and $q_k$ of faces $F_i$ and $F_k$, respectively, in a graph $\Gamma$.
By Lemma~\ref{triv} and Remark~\ref{mainremark}, $\Gamma$ can be transformed to a chain equivalent graph $\Gamma'$ with the same multiset of double perimeters, such that in the dual tree $T(\Gamma')$ the vertices of degrees $q_i$ and $q_k$ are adjacent. Without loss of generality, we are to increase $q_i$. In this case $q_i<P_i$ and $q_k \geq 3$. Consider faces $F_1$ and $F_2$ of $\Gamma'$ which can be described in terms of dual tree $T(\Gamma')$ as follows: $F_1=\pi_{F_k}^{-1}F_i$, $F_2=\pi^{-1}_{F_k}F_1$ (all faces $F_1$, $F_2$ and $F_i$ are different, since $q_k \geq 3$). Since $q_i<P_i$, there is at least one leaf in $F_i$.The following sequence of transformations finishes the proof (see Figure~\ref{fig:klever}; in the picture the shifts of leaves are omitted).
\begin{figure}[ht]
\begin{center}
		\includegraphics{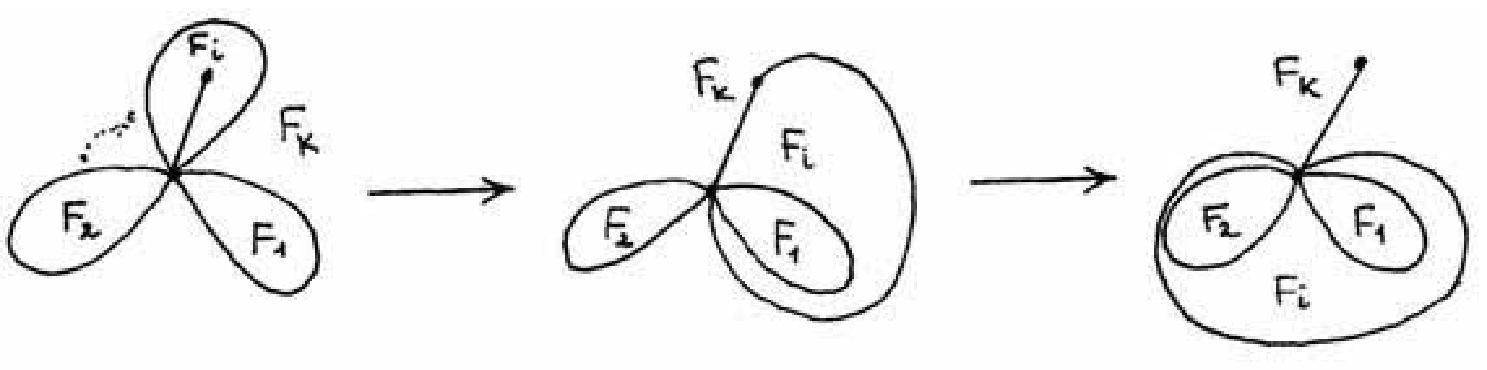}
		\caption{to Proposition~\ref{double}}
  	\label{fig:klever}
\end{center}
\end{figure}
Thus $q_i$ is increased by 2 and $q_k$ is decreased by $2$, which was required. 
\end{proof}

\noindent Altogether, we get
\begin{proposition}\label{itog1}
Two plane graphs with the same multiset of perimeters of faces are chain equivalent. 
\end{proposition}
Now we again consider graphs with labeled vertices, i.e., we return the multiplicities of vertices into consideration.
In statements from Lemma~\ref{lem1} to Theorem~\ref{maintheorem} we prove that {\it if two plane Brauer graphs with the same multisets of labels of vertices are isomorphic as non-labeled graphs, then they are chain equivalent as labeled graphs. }
In view of the above arguments, it's enough to prove this for reduced graphs. Moreover, in the case of bipartite graphs we may restrict ourselves to considering graphs in canonical form. Recall that the process of putting a graph to reduced form (and to canonical form, for graphs of type 1) started with choosing an {\it arbitrary} vertex $A$. Recall also that we can arbitrarily shift leaves in a face, by tilting transformation of type 1.

\begin{figure}[ht]
\begin{center}
		\includegraphics{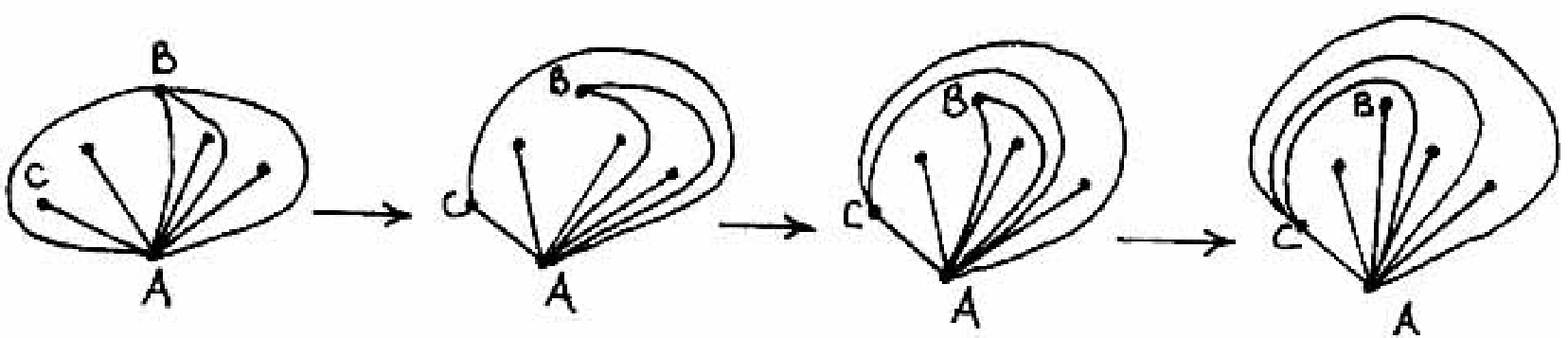}
		\caption{Proof of Lemma~\ref{lem1}}	
		\label{fig:lemma1}
	\end{center}
\end{figure}

\begin{figure}[ht]
\begin{center}
		\includegraphics{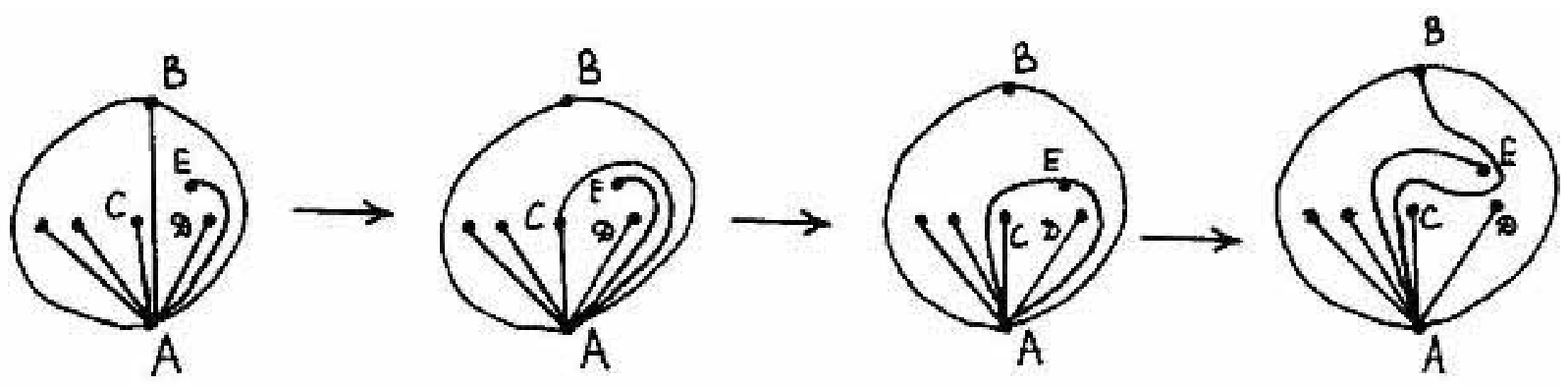}
		\caption{Proof of Lemma~\ref{lem2}}
		\label{fig:lemma2}
	\end{center}
\end{figure}

\noindent{\bf I. Reduced graphs of type 1.} 
For reduced graphs of type 1, it suffices to prove the following lemmas:
\begin{lemma} \label{lem1}
Let $\Gamma$ be a graph in canonical form, let $B\neq A$ be the second non-dangling vertex of $\Gamma$, let $F$ be a face. Then $\Gamma$ is chain equivalent to a graph in canonical form, in which
\begin{enumerate}
	\item $B$ is a dangling vertex in the face $F$.
	\item Some vertex $C \neq A$ which belongs in $\Gamma$ to $F$ is a non-dangling vertex.
	\item The other dangling vertices belong in $\Gamma$ and in $\Gamma_1$ to the same faces. 
\end{enumerate}
\end{lemma}
\begin{lemma} \label{lem2}
Let $\Gamma$ be a graph in canonical form, let $B\neq A$ be the second non-dangling vertex of $\Gamma$, let faces $F_1$ and $F_2$ be adjacent. Then $\Gamma$ is chain equivalent to a graph in canonical form $\Gamma_1$, in which 
\begin{enumerate}
	\item There is a dangling vertex which belongs in $\Gamma$ to $F_1$ and belongs in $\Gamma_1$ to $F_1$, and there is another dangling vertex which belongs in $\Gamma$ to $F_2$ and in $\Gamma_1$ to $F_2$. 
	\item The other dangling vertices belong in $\Gamma$ an in $\Gamma_1$ to the same faces. 
\end{enumerate}
\end{lemma}
\noindent For the proof of Lemma~\ref{lem1} see Figure~\ref{fig:lemma1}. For the proof of Lemma~\ref{lem2} see Figure~\ref{fig:lemma2}.

\begin{figure}[ht]
\begin{center}
		\includegraphics{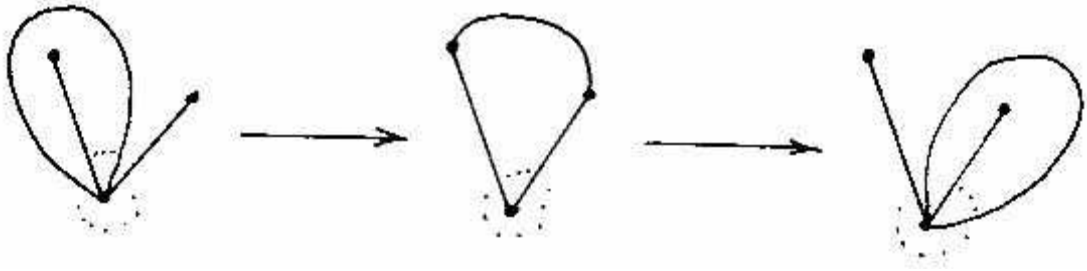}
		\caption{Interchange between adjacent faces}
  	\label{fig:short}
\end{center}
\end{figure}

\noindent{\bf II. Reduced graphs of type 2.}
Since the dangling vertices in a face can be shifted in arbitrary way, it's enough to show how to interchange dangling vertices belonging to different faces (say, to $F_1$ and $F_2$). First consider the case when the external perimeter of $F_1$ or $F_2$ is greater then 1. Then by Lemma~\ref{triv} and Remark~\ref{mainremark}, there is a sequence of tilting transformations making $F_1$ and $F_2$ adjacent. Moreover, this sequence preserves the faces to which belong the dangling vertices (see Figure~\ref{fig:butterfly}). Therefore, in this case it's enough to show how to interchange dangling vertices which belong to adjacent faces: see Figure~\ref{fig:short}.
\noindent Now consider the case when the dangling vertices which we want to interchange belong to faces, which correspond to dangling vertices of the dual tree. 
\begin{lemma}
Let $T$ be a plane tree with labeled vertices, let $V_1$ and $V_2$ be dangling vertices of $T$. Suppose that $T$ is not a chain. Then $T$ is equivalent to a tree, in which the edges which are incident with $V_1$ and $V_2$ are incident to a common vertex $V$. Moreover, $\pi_{V}(VV_1)=VV_2$.
\end{lemma}
\begin{proof}
By Remark~\ref{mainremark} it is enough to find a tree with the same multiset of degrees as in $T$, in which some two leaves are adjacent to a common vertex. Denote the degrees of $T$ by $r_1,\dots, r_g$ in such way that $r_1=r_2=1$, $r_3\geq 3$. Observe that the sum of numbers $d_3-2,d_4,\dots,d_g$ equals $2g-6$. It can be shown by induction on $g$ that there is a tree $T'$, in which these numbers are the degrees of vertices. To obtain the needed tree, we add two leaves to the vertex of $T'$ of degree $r_3-2$.  
\end{proof}
\begin{figure}[!ht]
\begin{center}
		\includegraphics[scale=0.9]{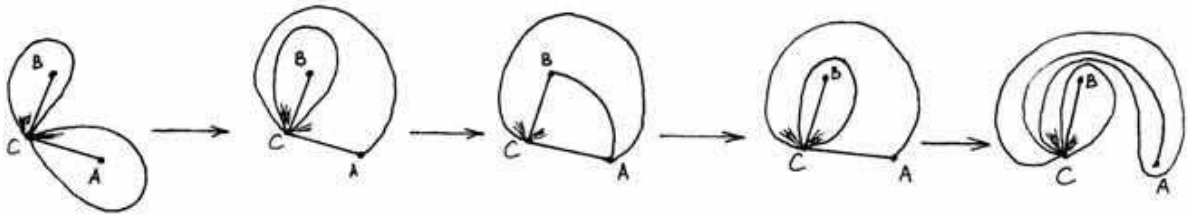}
		\caption{Interchange between "dangling" faces}
  	\label{fig:long1}
\end{center}
\end{figure} 
We see that if $T$ is not a chain, then it suffices to show how to interchange dangling vertices between two "dangling" faces, which have a common adjacent face: see Figure~\ref{fig:long1}.
\begin{figure}[!ht]
\begin{center}
		\includegraphics{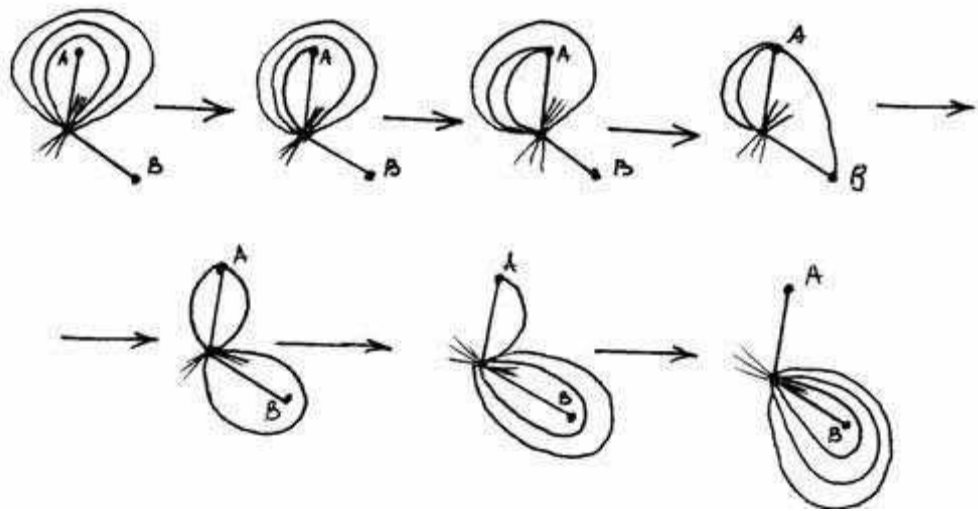}
		\caption{Chain case}
  	\label{fig:long2}
\end{center}
\end{figure} 
It remains to examine the case when $T(\Gamma)$ is a chain, and $F_1$ and $F_2$ correspond to the two dangling vertices of $T(\Gamma)$. If some other face of $\Gamma$ contains a dangling vertex, the needed interchange comes to three interchanges of the above form. In Figure~\ref{fig:long2} is is shown how to interchange leafs in case when the rest faces have perimeter 2. (For the graph in the picture $g=4$, and this case fully represents the general case.) 

\noindent This finishes the proof of the main theorem in this section:

\begin{theorem} \label{maintheorem}
Let $\Lambda_1$ and $\Lambda_2$ be symmetric $SB$-algebras of genus $0$. Then $\Lambda_1$ and $\Lambda_2$ are derived equivalent if and only it their Brauer complexes have the same multisets of perimeters of faces and the same multisets of labels on vertices.
\end{theorem}


\end{document}